\documentclass[12pt, A4paper]{article}
\usepackage[utf8]{inputenc}
\usepackage{amsthm,amsmath,amssymb,graphicx,enumerate,comment}

\usepackage{xcolor,colortbl}

\newtheorem{theorem}{Theorem}[section]
\newtheorem{lemma}[theorem]{Lemma}
\newtheorem{proposition}[theorem]{Proposition}
\newtheorem{corollary}[theorem]{Corollary}

\newtheorem{conjecture}{Conjecture}

\theoremstyle{definition}

\newtheorem{remark}{Remark}

\newcommand{\ZZ}{\mathbb{Z}}
\newcommand{\B}{{\mathcal B}}

\newcommand{\arxiv}[1]{\texttt{arXiv:#1}}

\title{On the Hamiltonian Bicirculants}

\author{Simona\,Bonvicini, Toma\v z\, Pisanski and Arjana\,\v Zitnik }

\date{October 27, 2025}

\begin{document}

\maketitle

\begin{abstract}

A bicirculant is a regular  graph  that admits a semi-regular automorphism  with two vertex-orbits of the same size. By  $m$ we denote  the size of vertex-orbits and by $d$ the valence of a bicirculant. Furthermore, we denote by $s$  the valence of the bipartite graph joining the two vertex-orbits. In 1983, Brian Alspach proved that the only non-hamiltonian generalized Petersen graphs are $G(m,2)$ with $m \equiv 5 \pmod 6$. In a recent paper we conjectured  that this is the only exception among regular, connected bicirculants of degree $d > 1$ and we have verified the conjecture for the quartic bicirculants with $s=2$, also known as the generalized rose window graphs.

In this paper we develop tools and apply them for a partial verification of the conjecture. We show that the conjecture holds for all bicirculants with $s \leq 2$. As a consequence we obtain that every connected bicirculant with $s \ge 3$ is hamiltonian if $m$ is  a product of at most three
prime powers.  In particular, every connected bicirculant with $s \ge 3$ is hamiltonian for even $m<210$ and odd $m < 1155$. 
Our results imply that many other families of bicirculants are hamiltonian. For example, all bicirculants with $d-s$  odd are hamiltonian.
\end{abstract}


\textbf{Math Subj Class (2020)}: 
05C45, %Eulerian and Hamiltonian graphs
05C25, %Graphs and abstract algebra (groups, rings, fields, etc.) 
05C76,  %Graph operations (line graphs, products, etc.)
05C70,  % 	Edge subsets with special properties (factorization, matching, partitioning, covering and packing, etc.)
05E18. % group actions on combinatorial structures
%20B25. %finite automorphism groups of algebraic, geometric or combinatorial structures

%%%%%%%%%%%%%%%%%%%%%%%%%%%%%%%%%%%%%%%%%%%%%%%%%%%%%%%%%%%%%%%%%%%%

\section{Introduction}

A \emph{bicirculant} is a regular graph that is a cyclic cover over a two-vertex (pre)-graph. In an equivalent way, bicirculants are defined to be regular graphs that admit an automorphism having two vertex-orbits of the same size \cite{Pi2007}.  
A typical bicirculant can be described as follows; see for example \cite{MALNIC2007891}. Given an integer $m \ge 1$ and sets $R,S,T \subseteq \ZZ_m$
such that $R=-R$, $T=-T$,   $0 \not\in R \cup T$   and $0 \in S$,  the graph $B(m;R,S,T)$ has vertex set 
$V=V_1 \cup V_2$, where $V_1=\{u_0,\dots,u_{m-1}\}$ and $V_2=\{v_0,\dots,v_{m-1}\}$, and edge set 
$$E=\{u_iu_{i+j}| \ i \in\ZZ_m, j \in R\} \cup 
      \{v_iv_{i+j}| \ i \in\ZZ_m, j \in T\} \cup
      \{u_iv_{i+j}| \ i \in\ZZ_m, j \in S\}.$$ 
Obviously, the mapping $\alpha:V \to V$, defined by $\alpha(u_i)=u_{i+1}$, $\alpha(v_i)=v_{i+1}$
is an automorphism of $B(m;R,S,T)$, having two vertex-orbits of the same size.

We call the vertices from $V_1$ the \emph{outer vertices} and the vertices from $V_2$ the \emph{inner vertices}.
There are three types of edges: the edges incident to two outer vertices are called \emph{outer edges},
the edges incident to two inner vertices are called \emph{inner edges}, and edges connecting an outer vertex to an inner vertex are called  \emph{spokes}. Specifically, the edges $u_iu_{i+a}$, $i \in \ZZ_m$, $a \in R$, are called \emph{outer edges of type} $a$, the edges $v_iv_{i+b}$, $i \in \ZZ_m$, $b \in T$, are called \emph{inner edges of type} $b$ and the edges $u_iv_{i+c}$, $i \in \ZZ_m$, $c \in S$, are called \emph{spokes of type} $c$. Furthermore, we call the subgraph of $B(m;R,S,T)$ induced by the outer (inner) vertices the \emph{outer rim} (the \emph{inner rim}, respectively). The inner rim and the outer rim are circulants, described by the sets $R$  and $T$, respectively. The bipartite subgraph of $B(m;R,S,T)$, spanned by the spokes, is  a \emph{cyclic Haar graph}, which we will denote by $H(m;S)$; so $H(m;S)= B(m;\emptyset,S,\emptyset)$. 
The family of cyclic Haar graphs was introduced by Hladnik, Maru\v si\v c and Pisanski in  2002 \cite{HlMaPi2002}.

We consider only regular bicirculants, so $|R|=|T|$, and we denote $|R| = |T| = r$, $|S| = s$. The degree of such a bicirculant is then $d=r+s$. We denote by $\B(m;d,s)$  the class of $d$-valent bicirculants of order $2m$ with spoke edges forming a bipartite subgraph of degree $s$. 

\begin{figure}[htb]
    \centering
    \includegraphics[scale=0.7,angle=0]{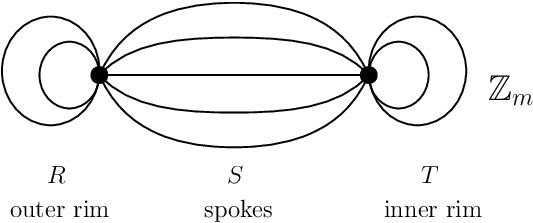}
    \quad
    \includegraphics[scale=0.7,angle=0]{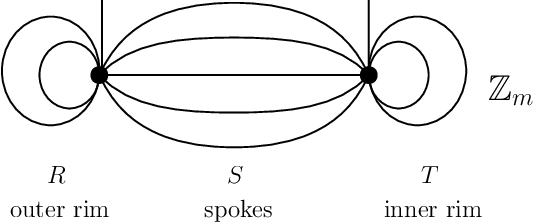}
    \caption{A typical bicirculant  $B(m;R,S,T) \in \B(m;d,s)$, defined by a voltage graph for $d-s$ even (left); for $d-s$ odd a semi-edge is drawn at each vertex (right). }
    \label{fig:my_label2}
\end{figure}

Bicirculants  can also be described as  regular $\ZZ_m$-covers  over   graphs on two vertices with possible multiple edges, loops and semi-edges; see \cite{Pi2007}. Figure  \ref{fig:my_label2} shows a typical regular voltage graph on two vertices that defines a bicirculant.  We see that there are two cases. On the left, the voltage graph has no semi-edges. On the right, the voltage graph has two semi-edges, one at each vertex. The latter case occurs if and only if $m$ is even  and $r = d -s$ is odd. In this case, the voltages on the semi-edges are equal to $m/2$. 

Many well-known graphs are bicirculants. For example, 
the generalized Petersen graphs are bicirculants, of valency 3. More generally, an $I$-graph $I(m;j,k)$ is a bicirculant $B(m;R,S,T)$ with $m \ge 3$, $R=\{j,-j\}$, $T=\{k,-k\}$ and $S=\{0\}$, where $0 < j,k < m/2$; see  \cite{Foster}. 
The $I$-graphs thus belong to the class $\B(m;3,1)$.
A generalized rose window graph  \cite{BoPiZi2025,wilsonRW}  is a bicirculant graph of degree 4  from the class $\B(m;4,2)$, obtained from an $I$-graph by adding an additional set of spokes.
The  pentavalent generalized  Taba\v cjn graphs \cite{Taba1,Taba2} form the class $\B(m;5,3)$.  Also Cayley graphs on dihedral groups are bicirculants $B(m;R,S,T)$ with $R=T$.

\medskip

Several papers investigated various properties of certain subfamilies of bicirculant graphs, such as generalized Petersen graphs, rose window graphs, Taba\v cjn graphs, cyclic Haar graphs, etc. Our aim is different. We are considering a classical problem of the existence of a hamilton cycle, which is known to be NP-complete, for general graphs. However, we hope that an effective classification of hamiltonian bicirculants is possible. Our approach is very optimistic. Namely, we believe that all non-hamiltonian bicirculants are already known and hope that for all other bicirculants efficient constructions of hamilton cycles can be found. Our strategy is based on a useful property of bicirculants. Namely, by removing a suitable $1$-factor or $2$-factor, a $d$-valent  bicirculant reduces to a $(d-1)$-valent or, respectively, to a $(d-2)$-valent bicirculant. The latter may be disconnected. However, even in such a case its connected components are bicirculants. The main idea is to apply various constructions that link hamilton paths in each component to a hamilton cycle of the original graph.

In a recent paper \cite{BoPiZi2025} we proved that all generalized rose window graphs are hamiltonian. Based on this fact and the fact that the classification of Brian Alspach of hamiltonian generalized Petersen graphs \cite{Al11983}  has been extended by Bonvicini and Pisanski to $I$-graphs \cite{BoPi2017}  and that no new 
non-hamiltonian graphs have been found, we conjectured that the only non-hamiltonian bicirculants are the $K_2$ and the generalized Petersen graphs, determined by Alspach.

\begin{conjecture} \label{conject}
Every connected bicirculant, except for the complete graph $K_2$ and the generalized Petersen graphs $G(m, 2)$ with $m\equiv 5\pmod 6$,  is hamiltonian. \end{conjecture}

In this paper we give a partial solution to Conjecture  \ref{conject}. 
After reviewing known results on hamiltonicity and basic properties of bicirculants  (Section \ref{section:properties}), we  show how to combine hamilton cycles in the connected components of certain subgraphs of a connected bicirculant to a hamilton cycle of the whole graph (Section \ref{sec:hamconstrucion}). 
We use constructions from Section \ref{sec:hamconstrucion} to prove that certain families of bicirculants are hamiltonian.
Our main results  are summarized below and they are proved in Sections \ref{subsection:consequences} and \ref{section:main}.

\begin{theorem}\label{th_bicirculant_s2}
Every connected bicirculant graph whose vertices are incident to $s\le 2$ spokes is hamiltonian, except for 
$K_2$ and  the generalized Petersen graphs $G(m, 2)$ with $m\equiv 5\pmod 6$.
\end{theorem}

\begin{theorem}\label{thm:primes}
Let $G=B(m;R,S,T)$ be a connected bicirculant  with  $|S| \ge 3$. If $m$ is a product of at most three prime powers, then $G$ is hamiltonian.

In particular, every connected bicirculant graph of order $2m$ is hamiltonian for  even  $m<210$ %$m<2\,310$ 
and odd $m< 1155$, with the exception of the complete graph $K_2$ and  the generalized Petersen graphs $G(m, 2)$ with $m\equiv 5\pmod 6$.
\end{theorem}

\begin{theorem}\label{th_m/2}
Let  $G=B(m; R,  S, T)$ be a connected bicirculant graph with $m$ even, $m\ge 4$, $|S|\ge 3$ and $m/2\in R, T$. Then  $G$ is hamiltonian.
\end{theorem}

\begin{theorem}\label{th_main}
Let $G=B(m; R, S, T)$ be a connected bicirculant graph  with $m\ge 3$ and $|S|\ge 3$.
Assume that the connected components of $H(m; S)$ are hamiltonian and $R\cup T$ contains at least one element that is not coprime to
$\gcd(m, S)$ in the case where $\gcd(m, S)>1$. Then $G$ is hamiltonian.
\end{theorem}

\section{Background and preliminary results}
\label{section:properties}

In this section we summarize some known results on the hamiltonicity and basic properties of bicirculant graphs that are relevant for our paper. We note that  bicirculants $B(m;R,S,T)$ with $R=T$ are Cayley graphs on dihedral groups. In particular, cyclic Haar graphs are also Cayley graphs on dihedral groups.

In \cite{Al11983} Brian Alspach completely resolved the problem which generalized Petersen graphs are hamiltonian. Later   it was shown by Bonvicini and Pisanski \cite{BoPi2017} that the nonhamiltonian generalized Petersen graphs are the only nonhamiltonian $I$-graphs.

Alspach and Zhang \cite[Theorem 3.1]{AlZh1989} have shown that all cubic Cayley graphs over dihedral groups are hamiltonian. From their result it follows that all connected cubic cyclic Haar graphs are hamiltonian. 

\begin{theorem}\label{th_cyclic_Haar_s3}
Every connected cubic cyclic Haar graph is hamiltonian.
\end{theorem}

Consequently, we have the following result that completely settles the problem which cubic bicirculants are hamiltonian; see also \cite{BoPiZi2025}.

\begin{theorem} 
\label{th_cubicbicirculants}
Every connected cubic bicirculant graph is hamiltonian, except for the generalized Petersen graphs $G(m,2)$ with $ m \equiv 5 \pmod 6$.
\end{theorem}

Since the family of generalized Petersen graphs  $G(m,2), m \equiv 5  \pmod 6$,  was first discovered by Brian Alspach, we name these graphs \emph{the Alspach generalized Petersen graphs}.

Other families of bicirculants for which the hamiltonian property was studied include Cayley graphs on dihedral groups and the generalized rose window graphs. 
For convenience, we recall the result  of Alspach, Chen and Dean \cite{AlChDe2010}.

\begin{theorem}[\protect{\cite[Theorem 1.8]{AlChDe2010}}]\label{th_generalized_dihedral}
A Cayley graph of valency at least three on a generalized dihedral group, whose order is divisible by four, is hamilton-connected, unless it is bipartite, in which case it is hamilton-laceable.
\end{theorem}

We will also use the following theorem from our recent paper.

\begin{theorem}[\cite{BoPiZi2025}] \label{th_rwgraphs}
Every connected generalized rose window graph is hamiltonian.
\end{theorem}

Now we recall some properties of bicirculants.  
From the theory of covering graphs it follows that we may choose the elements of the sets $R,S,T$ to be such that $0 \in S$, which we have already assumed. In this case we have a simple criterion when  a bicirculant is connected; see, for example, \cite{GT}.

\begin{proposition}  \label{prop:connected}
A bicirculant graph  $B(m;R,S,T)$  is connected if and only if $\gcd(m,R,S,T)=1$. 
\end{proposition}

In the case that a bicirculant graph is disconnected, it is composed of isomorphic connected components. We will use the following notation. Let $A$ be a set and let $i$ be an integer. We define $A-i=\{a-i \ | \ a \in A \}$ and $A/i= \{a/i \ | \ a \in A \}$.

\begin{proposition}\label{pro:d_connected_components} 
Let $G = B(m;R,S,T)$. Suppose $\delta = \gcd(m,R,S,T)>1$. 
Then $G$ is a disjoint union of $\delta$   isomorphic graphs $G_0, \dots G_{\delta-1}$ such that $u_i \in G_i$, $i = 0,\dots,\delta-1$. Moreover, each graph $G_i$  is connected and isomorphic to the graph $B(m/\delta;R/\delta,S/\delta,T/\delta)$. 
\end{proposition}

Different parameters may give the same bicirculants. In particular, it is easy to verify that the following result holds.

\begin{lemma}  \label{lemma:iso}
Let $G = B(m;R,S,T)$ be a bicirculant and let $c \in S$. Then the graph $B(m;R,S-c,T)$ is isomorphic to the graph $G$.
\end{lemma}

Next  we consider relationship between pairs of bicirculants of order $2m$ 
that only differ in the presence of inner and outer edges of type $m/2$. Let $m$ be even and let $G = B(m; R, S, T)$ be a bicirculant such that $m/2  \not \in R \cup T$. We denote by $G^=$ the graph $B(m; R \cup \{m/2\}, S, T \cup \{m/2\})$.

\begin{lemma} \label{lemma_cprod}
Let $m$ be even and let $G = B(m; R, S, T)$ be a  bicirculant such that $m/2  \not \in R \cup T$. Suppose that $G^=$ is  connected.  Then either $G$ is connected, or $G$ is composed of two isomorphic  connected components $G_0$.
In the latter case, $G^=$ is isomorphic to the cartesian product 
$G_0 \square K_2$.
\end{lemma}

\begin{proof}
If $G$ is connected, there is nothing to prove, so assume that $G$ is disconnected. Then $\delta = \gcd(m,R,S,T) > 1$. Since $\gcd(m/2,R,S,T) = 1$ it follows that $\delta = 2$. Hence $G$ consists of two isomorphic copies of a connected graph $G_0$. 
The outer/inner edges of type $m/2$ determine a perfect matching in $G^=$ that connects vertices of one copy of $G_0$ to the corresponding vertices of the second copy. This means that $G^=$ is the cartesian product of $G_0$ and $K_2$.
\end{proof}

\begin{lemma} \label{lemma_cprodham}
Let $m$ be even and let $G = B(m; R, S, T)$ be a connected  bicirculant such that $m/2 \in R \cap T$. Let $R_0 = R \smallsetminus \{ m/2\}$, $T_0 = T \smallsetminus \{m/2\}$ and $H = B(m;R_0,S,T_0)$.

\begin{itemize}
\item[(i)]  The graph $H^=$ is isomorphic to the graph $G$. 
\item[(ii)] If $H$ is   hamiltonian, then also   $G$ is hamiltonian.
\item[(iii)] If $H$ is disconnected,  then it is composed of two isomorphic copies of a connected graph $H_0=B(m/2; R_0/2,S/2,T_0/2)$. If $H_0$ has a hamilton path, then $G$ is hamiltonian.
\end{itemize}
\end{lemma}
\begin{proof}
\noindent
(i) The claim follows from the definition of $H^=$.\\ 
(ii) The claim follows from   the fact that $H$ is a spanning subgraph of $G$.\\
(iii) By Lemma \ref{lemma_cprod}, graph $G$ is isomorphic to the cartesian product of $H_0$ and $K_2$. A cartesian product of a graph with hamilton path and $K_2$ is always hamiltonian, which is well known and easy to prove. 
\end{proof}

\section{Some hamiltonian constructions}
\label{sec:hamconstrucion}

In this section we develop some methods for the construction of a hamilton cycle in a connected bicirculant graph $B(m; R, S, T)$, with $R$, $T$ non-empty. 
For convenience we will denote a bicirculant $B(m;\{a,-a\},S,\{b,-b\})$ simply by $B(m;a,S,b)$.

The first method of construction is described in Lemma \ref{pro_method1_case1}.
It joins  hamilton  cycles in the connected components of a cyclic 
Haar graph $H(m; S)$ into a hamilton cycle of a connected bicirculant graph  $B(m; a, S, b)$ containing $H(m; S)$ as a disconnected spanning subgraph.  This construction can be used when   at least one of the integers $a,b$  is different from $m/2$ and  not coprime to $\gcd(m,S)$.
 
 The second method of construction is described in Lemma \ref{pro_construction_Bd1bis}:
 it joins   hamilton cycles in the connected components of $B(m; R\smallsetminus\{a,-a\},  S, T\smallsetminus\{b,-b\})$ into a hamilton cycle of
 $B(m; R, S, T)$, where $a \in R$, $b \in T$ and $a,b \ne m/2$.  Unlike the hamilton cycles used in Lemma \ref{pro_method1_case1},  the hamilton cycles employed in Lemma \ref{pro_construction_Bd1bis} contain outer and inner edges; they are joined through outer and inner edges of types $a$ and $b$ in both constructions.

 As a consequence of the combined application of Lemma \ref{pro_method1_case1}  and Lemma \ref{pro_construction_Bd1bis},
 we  show that a connected bicirculant graph whose vertices are incident to at most two spokes is hamiltonian, except for the Alspach generalized Petersen graphs and the complete graph $K_2$;
 see Section \ref{subsection:consequences}.
 
The construction of a hamilton cycle in $B(m; a, S, b)$ given  in Lemma \ref{pro_method1_case1} is an alternative to the construction given for the generalized rose window graphs in \cite{BoPiZi2025},  but it does not allow to find  a hamilton cycle in every generalized rose window graph because it can be applied only when at least one of the integers $a$, $b$ is not coprime to $\gcd(m, S)$. 

\smallskip

To simplify the proof of  the constructions, we first give a result on bicirculant graphs of order $\le 10$.

\begin{lemma}\label{lemma_bicirculant_m<5}
Every connected bicirculant graph  of order $2m$ with $m\le 5$ is hamiltonian, except for the  complete graph $K_2$ and the 
Petersen graph $G(5,2)$. 
\end{lemma}

\begin{proof}
Let $G=B(m; R, S, T)$ be a connected bicirculant graph of order $2m$ with $ m\le 5$.
If $m=1$, then $G$ is the complete graph $K_2$, which is not hamiltonian.
If $m=2$,   then   $G$ contains a $4$-cycle, hence it is hamiltonian.  
Consider the case when  $m \in \{3,5\}$ and assume that the graph $G$ is not isomorphic to the Petersen graph. If $|R| \ge 1$, then the graph $G$ contains a subgraph isomorphic to the prism $B(m;1,\{0\},1)$  
which is hamiltonian.  Otherwise the set $S$ contains at least two elements, say $0$ and $c \in \{1,\dots,m-1\}$, and the subgraph of $G$ induced by the spokes of types $0$ and $c$ is a hamilton cycle of $G$ since $\gcd(m, c)=1$. 

Now consider the case $m=4$.
If $1 \in R$, then also $1 \in T$ and the graph $G$ contains  a subgraph isomorphic to the prism $B(4;1,\{0\},1)$, which is hamiltonian.
Otherwise $\{0,c\} \subseteq S$ for some $c \in \{1,3\} \subseteq S$, since the graph $G$ is connected.  Now the subgraph of $G$ induced by the spokes of type $0$ and $c$ is a hamilton cycle of $G$.
\end{proof}

\subsection{Construction for the graphs $B(m; a, S, b)$}
\label{subsection:construction1}

In order to apply our first construction,  we will represent a bicirculant graph as we are going to describe. Let $G=B(m; a, S, b)$ be a connected bicirculant graph with $m \ge 2$, $|S|\geq 2$, $\gcd(m, S)>1$,   $a,b\ne m/2$,  and at least one of the integers $a$, $b$, say $b$,  is not coprime to $\gcd(m,S)$. 

\smallskip

We denote with $H=H(m; S)$  the cyclic Haar graph induced by the spokes of $G$, 
and with  $K$ the subgraph obtained from $G$ by removing the edges of type $a$ from $G$, that is,  $K$ is the  graph $B(m; \emptyset,S,   \{b,-b\})$. The connected components of $K$ form a partition of the components of $H$: each component of $K$ consists of $\gcd(m, S)/\gcd(m, S, b)$ components of $H$, which are connected by edges of type $b$.  We set 
$$\lambda=\gcd(m, S, b)-1\mbox{ and } \quad\mu=\gcd(m, S)/\gcd(m, S, b)-1$$
and denote with $K_0,\ldots, K_{\lambda}$ the connected components of $K$; by the assumptions we have $\lambda >0$.
For $0\le j\le\lambda$, we denote by $H_{i,j}$, with $0\le i\le\mu$ , the connected components of $H$ that are contained in $K_j$.
Moreover,  $H_{0,0}$ will be the component of $H$ containing the vertex $u_0$, and consequently $K_0$ will be the component of $K$ containing $u_0$. 
In the following we will use the notation $x\,P\,y$ to denote  a subpath of a path $P$, starting at  vertex $x$ and ending at  vertex $y$.

We will assume that the connected components of $H$ are hamiltonian, therefore they also contain  hamilton paths. Let $P$ be a hamilton path in $H_{0,0}$.
We can assume that $P$ is a hamilton path with origin in $v_s$ and terminus in $u_0$; we denote with $u_z$ an arbitrary vertex of $P$ different from $u_0$
and with $v_h$, $v_k$ the vertices adjacent to $u_z$ in $P$.  For convenience of notation, the remaining vertices of $P$ will be simply denoted with $u_i$, $v_i$ instead of $u_{x_i}$, $v_{y_i}$.
More specifically, we set 
$$P=(v_s, u_t, v_t, u_{t-1}, \ldots,   v_h,  u_z,  v_k,  \ldots,   u_1,  v_0 , u_0).$$
The paths $v_s\,P\,u_z$ and $u_z\,P\,u_0$ can have arbitrary length $>0$. It may happen that $v_0=v_k$ and $u_1=u_z$ and/or $v_s=v_h$ and $u_t=u_z$.
It is understood that $u_0$ and $v_s$ are adjacent if $P$ provides a hamilton cycle in $H_{0,0}$,  
which is the case that will be considered in the construction we are going to define.

The connected components of $K$ -- and also of $H$ -- contain copies of $P$ that, for our purposes, we will describe as follows.
For $0\leq i\leq\mu$ and $0\leq j\leq\lambda$, we denote with $P_{i,j}$ the path obtained from $P$ by adding $ib+ja\pmod m$ to the subscripts of the vertices in $P$;
$P_{i,j}$ is the copy of the path $P$ in the component $H_{i,j}$ of $H$, and $P_{0,0}$ corresponds to the path $P$. 
Accordingly, the vertices in $P_{i,j}$ will be denoted with $u^{i,j}_x$, $v^{i,j}_x$, where $u^{i,j}_x=u_{x+ib+ja}$, $v^{i,j}_x=v_{x+ib+ja}$.  The outer vertices $u^{i,j}_x$ in $P_{i,j}$ are adjacent to the outer vertices $u^{i,j+1}_x$ in $P_{i, j+1}$ by the edges of type $a$; the inner vertices
$v^{i,j}_x$ in $P_{i,j}$ are adjacent to the inner vertices $v^{i+1,j}_x$ in $P_{i+1, j}$ by the edges of type $b$. 
For simplicity of notation, we will denote with $u_x\,P_{i,j}\,v_y$, or with  $u_x\,P\,v_y$ in $H_{i,j}$, the path corresponding to the subpath $u_x\,P\,v_y$ from the vertex $u_x$ to the vertex $v_y$ of $P$; 
it is understood that the vertices $u_x$, $v_y$ in the notation $u_x\,P_{i,j}\,v_y$, or in the notation $u_x\,P\,v_y$ in $H_{i,j}$,  correspond to the vertices $u^{i,j}_x$, $v^{i,j}_y$.  

In this setting we can represent each component $K_j$, with $0\leq j\leq \lambda$,  as the graph induced by the union of the paths $P_{i,j}$, with $0\leq i\leq\mu$.
Of course,  the  outer vertices of $K_j$ are connected to the outer vertices of $K_{j+1}$ through the edges of type $a$,  in order to form the whole graph $G$.
We can therefore think of $G$ as follows: arrange the  paths $P_{i,j}$ with the same index $i$ in the same $i$-th row, and the paths $P_{i,j}$ with the same index $j$
in the same $j$-th column. For this placement, the component $K_j$ will also be called `the vertical $j$-th component'. 

Now we introduce our methods of construction.  The key of the constructions presented in  Lemmas \ref{pro_method1_case1} and  \ref{pro_construction_Bd1bis} 
 is a set of prescribed paths that can be joined appropriately to obtain a hamilton cycle in the
original graph $G$.  The proofs of Lemmas \ref{pro_method1_case1} and  \ref{pro_construction_Bd1bis} 
describe in the detail the paths involved in the constructions and their connections; paths and connections are easily readable from the figures accompanying the proofs (Figures \ref{fig_path_lambda_positive_mu0}--\ref{fig_construction_Bd1}).

\begin{lemma}\label{pro_method1_case1}
Let $G=B(m; a, S, b)$ be a connected bicirculant graph with $m\ge 2$, $|S|\geq 2$, $\gcd(m, S)>1$,  $a,b\ne m/2$,  and $\lambda >0$.  
Assume that the connected components of $H(m; S)$ are hamiltonian.  Then $G$ is hamiltonian.
\end{lemma}

\begin{proof}
We set $C=P\cup u_0v_s$ to be a hamilton cycle of $H_{0,0}$, and distinguish the cases $\mu =0$ and $\mu >0$; in both cases 
we give a list of  subpaths of $P$ that are joined to form a hamilton cycle in $G$.   Moreover,  for $\mu >0$,  it is  sufficient to consider the case where $\lambda\equiv\mu\pmod 2$, and the case where $\lambda$ is even and $\mu$ is odd.  In fact,  due to the symmetry between outer and inner vertices, the case where  $\lambda$ is odd and $\mu$ is even can be obtained from the case where  $\lambda$ is even and $\mu$ is odd.

\medskip

\noindent
{\bf Case 1: $\mu=0$.}

First assume that $\lambda$ is even.  Each component $K_j$ of $K$ contains only the component $H_{0,j}$ of $H$, for $0\leq j\leq\lambda$; 
this means that the cycle $C=P\cup u_0v_s$  in $K_0 $ (or $H_{0,0}$) has chords of type $b$.  

We  denote with $u_r$, $u_{\ell}$ and $u_p$, $u_q$ the vertices adjacent to $v_{s-b}$ and $v_{s+b}$, respectively, and without loss of generality we can assume that they are arranged on $P$ in the order $ v_s$,  $u_t$,  $u_r$, $v_{s-b}$,  $u_{\ell}$,  $u_p$, $v_{s+b}$, $u_q$, $u_0$, (otherwise we swap $-b$ and $b$); see Figure \ref{fig_path_lambda_positive_mu0}(a). 
In this setting we can find a hamilton path in $H_{0,0}$ from $u_0$ to $u_p$ and a hamilton path in $H_{0,0}$ 
from $u_t$ to $u_{\ell}$, say $u_0\,P'\,u_p$ and $u_t\,P'\,u_{\ell}$.  More specifically,   $u_0\,P'\,u_p=u_0\,P\,v_{s+b}\cup\, v_{s+b}v_s\,\cup v_s\,P\,u_p$, and  
$u_t\,P'\,u_{\ell}=u_t\,P\,v_{s-b}\cup (v_{s-b}, v_s, u_0)\cup u_0\,P\,u_{\ell}$. 

We will also consider the path $u_0\,P''\,u_t$ from $u_0$ to $u_t$ given by $u_0\,P''\,u_t=u_0\,P\,v_{s+b}\cup (v_{s+b}, v_s,  v_{s-b})\cup v_{s-b}\,P\,u_t$; 
the union of the paths  $u_p\,P\,u_{\ell}$ and   $u_0\,P''\,u_t$  spans the vertices of $K_0$ (or $H_{0,0}$); a representation of the paths is given in Figure \ref{fig_path_lambda_positive_mu0}(a). 
It may happen that $u_t=u_r$ or $u_{\ell}=u_p$ or $u_q=u_0$. However, since $b \ne m/2$, the vertices $v_{s-b}$ and $v_{s+b}$ are distinct and they are also different from $v_s$. Therefore the graph $H_{0,0}$ contains at least   6 vertices that are needed in the definition of the above paths.

We use the above paths to construct a hamilton cycle in $G$.  More specifically, in $H_{0,0}$ we consider 
$u_0\,P'_{0,0}\,u_p$, in $H_{0,\lambda}$ we take $u_t\,P'_{0,\lambda}\,u_{\ell}$;
for $1\le j\le\lambda-1$,  we consider the paths $u_p\,P_{0,j}\,u_{\ell}$  and $u_0\,P''_{0,j}\,u_t$ in $H_{0,j}$. We find a hamilton cycle in $G$ by joining the paths through the following edges: the edges from $u_0$, $u_p$ in $H_{0,j}$ to $u_0$, $u_p$ in $H_{0,j+1}$, respectively, for $0\le j\le\lambda-1$, $j$ even, and the edges from $u_t$, $u_{\ell}$ in $H_{0,j}$ to $u_t$, $u_{\ell}$ in $H_{0,j+1}$, respectively, for $1\le j\le\lambda-1$, $j$ odd; Figure \ref{fig_path_lambda_positive_mu0}(b) shows the construction for $\lambda=4$. 

\begin{figure}[h]
\begin{center}
\includegraphics[width=13cm]{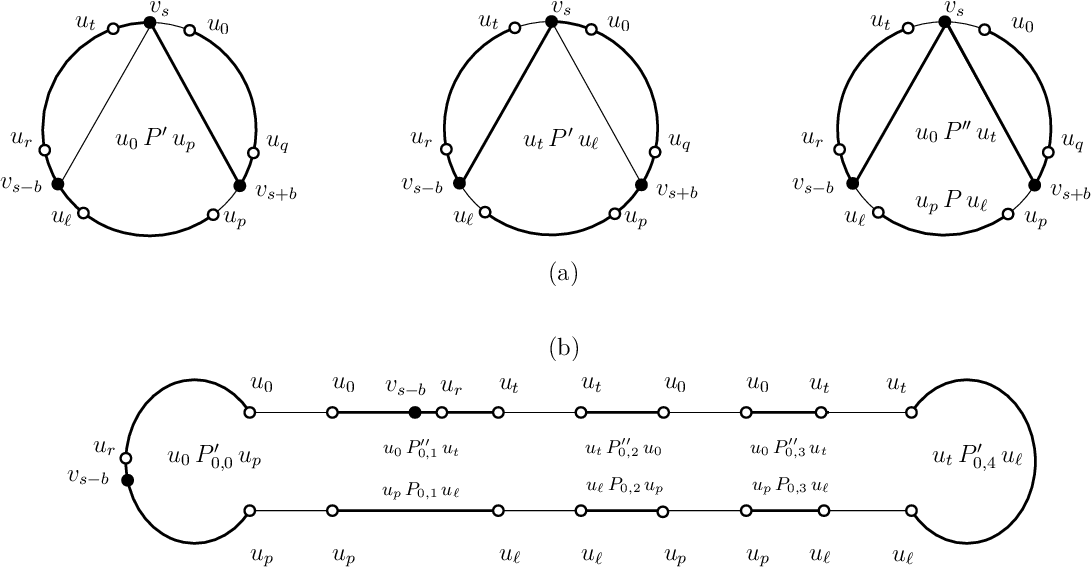}
\end{center}
\caption{The construction of a hamilton cycle in $G$ in Case 1 of the proof of Lemma  \ref{pro_method1_case1}, when $\lambda$ is even, $\mu=0$, and the connected components of $H(m,S)$ are hamiltonian.
(a) The bold lines stand for the paths we consider in each component $K_j$, $0\leq j\leq\lambda$. 
(b) We join the paths through edges of type $a$ -- represented by the thin lines -- as described in the proof of  Lemma  \ref{pro_method1_case1}; the figure shows the case $\lambda=4$. 
The same construction can be also used when  $\lambda$ is odd and $\mu=0$, by replacing the path $u_t\,P'_{0,\lambda}\,u_{\ell}$ with the path $u_0\,P'_{0,\lambda}\,u_p$.}\label{fig_path_lambda_positive_mu0}
\end{figure}

The previous construction also provides a hamilton cycle in $G$ 
when $\lambda$ is odd, $\lambda\ge 1$, and $\mu=0$: we simply replace the path $u_t\,P'_{0,\lambda}\,u_{\ell}$ in $H_{0,\lambda}$ with the path $u_0\,P'_{0,\lambda}\,u_p$;
for $\lambda=1$, it suffices to join the hamilton path $u_0\,P'\,u_p$ in $K_0$ to the hamilton path $u_0\,P'\,u_p$ in $K_{1}$ through the edges from $u_x$ in $K_0$ to $u_x$ in $K_1$, with $x\in \{0, p\}$. 

\medskip

\noindent
{\bf Case 2: $\mu>0$.}

The list of paths involved in the construction is the following:
\begin{itemize}
\item for $\lambda\ge 2$, in each vertical component $K_j$, with $1\leq j\leq\lambda-1$, we take the path $v_s\,P_{0,j}u_1$ and the edge $u^{0,j}_0v^{0,j}_0$;
for $1\leq i\leq\mu-1$,  we consider the paths $v_s\,P_{i,j}v_h$ and $v_k\,P_{i,j}v_0$ together with the isolated  vertices $u^{i,j}_0$ and $u^{i,j}_z$; 
if $\mu$ is even,  we take the isolated vertex $u^{\mu, j}_z$ and the path  $v_h\,C_{\mu,j}v_k=C_{\mu,j}-u^{\mu,j}_z$,  which is
the copy of the path $C-u_z$ that is obtained from $C=P\cup\,u_0v_s$ by deleting the vertex $u_z$;
if $\mu$ is odd,  we take the isolated vertex $u^{\mu, j}_0$ and the path  $v_s\,P_{\mu,j}v_0$;

\item in the vertical component $K_0$,  we will use the path $P_{0,0}$ from $v_s$ to $u_0$;
for $1\leq i\leq\mu-1$,  we take the path $v_s\,P_{i,0}v_h$ and $u_0\,P_{i,0}u_z$;  
if $\mu$ is even, we consider the  path $v_h\,C_{\mu,0}\,u_z=C_{\mu,0}-u^{\mu,0}_zv^{\mu,0}_h$,  which is the copy of the path $C-u_zv_h$ that is
obtained from $C=P\cup\,u_0v_s$ by deleting the edge $u_zv_h$;
if $\mu$ is odd,  we consider the path $P_{\mu,0}$, i.e., the copy of the path $P$ from $v_s$ to $u_0$ in $H_{\mu,0}$;

\item in the vertical component $K_{\lambda}$,  for $1\leq i\leq\mu-1$, we consider the paths $u_0\,C_{i,\lambda}v_h$ and $v_0\,C_{i,\lambda}u_z$, which are 
the copies of the paths that are obtained from $C=P\cup\, u_0v_s$ by removing the edges $u_0v_0$ and $u_zv_h$;  if $\lambda$ and $\mu$ are even, 
we take the path $v_h\,C_{\mu,\lambda}u_z$, and the path $u_1\,C_{0,\lambda}v_0=C_{0,\lambda}-u^{0,\lambda}_1v^{0,\lambda}_0$, 
 which  is the copy of the path $C-u_1v_0$ that is obtained from $C$ by removing the edge $u_1v_0$; 
 if $\lambda$ and $\mu$ are  odd,  we take the path  $u_0\,C_{\mu,\lambda}v_0$ and $u_0\,C_{0,\lambda}v_0$, which are the copies of the path
 $C-u_0v_0$ that is obtained from $C$ by deleting the edge $u_0v_0$; if $\lambda$ is even and $\mu$ is odd,  we take the path 
 $u_1\,C_{0,\lambda}v_0$ and the path  $u_0\,C_{\mu,\lambda}v_0$.
 \end{itemize}

Now we join the above paths to prove that $G$ is hamiltonian; the reader can refer to Figure \ref{fig_ham_cycle_lambda_mu_even_H(m,S)_hamiltonian} showing the case $\lambda=\mu=2$. For $\lambda\ge 2$, we consider the vertical component $K_j$,  with $1\leq j\leq\lambda-1$. For $0\leq i\leq\mu-1$, $i$ even,
we add the edges from the vertices $v_s$, $v_0$ in $H_{i,j}$ to the vertices $v_s$, $v_0$ in $H_{i+1, j}$, respectively; 
for $i$ odd,  we consider the edges from the vertices $v_h$, $v_k$ in $H_{i,j}$ to the vertices $v_h$, $v_k$ in $H_{i+1, j}$, respectively. In this way we find a path from $u_0$ in $H_{0,j}$ to $u_1$ in $H_{0,j}$ that covers all vertices of $K_j$,
with the exception of the vertices $u_0$,  $u_z$ in $H_{i,j}$ with $1\leq i\leq\mu-1$, and the vertex $u_{t'}$ in $H_{\mu,j}$,  where $u_{t'} = u_0$ if $\mu$ is odd,  $u_{t'} = u_z$ if $\mu$ is
even.

We denote such a path with $u_0\,K_j\,u_1$, and the set of uncovered vertices with $U_j$.  

\begin{figure}[h]
\begin{center}
\includegraphics[width=9cm]{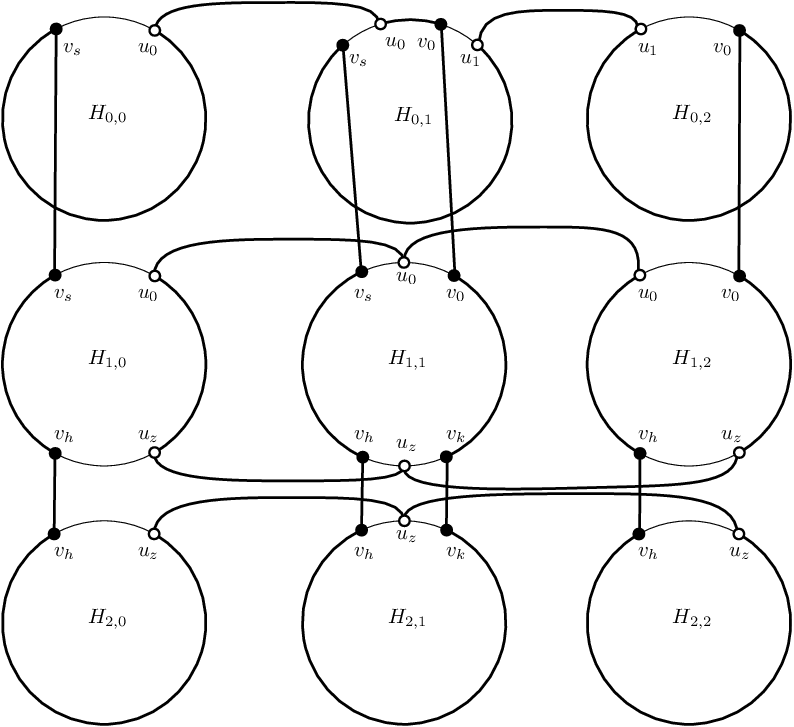}
\end{center}
\caption{The construction of a hamilton cycle in $G$ in Case 2 of the proof of Lemma  \ref{pro_method1_case1} when $\mu>0$; in the figure we have $\lambda=\mu=2$. }\label{fig_ham_cycle_lambda_mu_even_H(m,S)_hamiltonian}
\end{figure}

If $\lambda\ge 3$, we  connect $u_0\,K_j\,u_1$ to $u_0\,K_{j+1}\,u_1$ through the edges with both end-vertices equal to $u_0$ if $j$ is even, and through the edges with both end-vertices equal to $u_1$ if $j$ is odd, for $1\leq j\leq\lambda-2$.  

We thus construct the path $u^{0,1}_0\,P^*\,u^{0,\lambda-1}_{\alpha}$ from $u_0$ in  $H_{0,1}$  to $u_{\alpha}$ in $H_{0,\lambda-1}$, where $u_{\alpha} = u_0$ for $\lambda$ odd,  $u_{\alpha} = u_1$ for $\lambda$ even.  The path $u^{0,1}_0\,P^*\,u^{0,\lambda-1}_{\alpha}$
covers all vertices in $\cup^{\lambda-1}_{j=1} K_j$,  with the exception of the vertices in $\cup^{\lambda-1}_{j=1}\,U_j$.   
Note that the path $u^{0,1}_0\,P^*\,u^{0,\lambda-1}_{\alpha}$ is the path $u_0\,K_1\,u_1$ if $\lambda=2$, while there is no need to define it for $\lambda=1$, since a hamilton cycle in $G$ will be obtained from the paths in $K_0$ and $K_{\lambda}=K_1$, as we are going to explain.

We now construct a path covering all the vertices in $K_0\cup K_{\lambda}$, together with the vertices in $\cup^{\lambda-1}_{j=1}\,U_j$. 
We connect the subpaths defined in $K_0$ and $K_{\lambda}$ through the following edges: for $0\leq i\leq\mu-1$, $i$ even, we add the edges from $v_s$ in $H_{i,0}$ to $v_s$ in $H_{i+1, 0}$, together with the edges from $v_0$ in $H_{i,\lambda}$ to $v_0$ in $H_{i+1, \lambda}$; for $i$ odd,  we add the edges from $v_h$ in $H_{i,0}$ to $v_h$ in $H_{i+1, 0}$, together with
the edges from $v_h$ in $H_{i,\lambda}$ to $v_h$ in $H_{i+1, \lambda}$.  Further to the above,  for $1\leq i\leq\mu-1$ and $0\leq j\leq\lambda-1$,  we also add the edges from $u_0$, $u_z$ in $H_{i, j}$ to $u_0$, $u_z$ in $H_{i, j+1}$;
for $\mu$ odd, we add the edges from $u_0$ in $H_{\mu, j}$ to $u_0$ in $H_{\mu, j+1}$, with $0\leq j\leq\lambda-1$;
for $\mu$ even, we add the  edges from $u_z$ in $H_{\mu, j}$ to $u_z$ in $H_{\mu, j+1}$, with $0\leq j\leq\lambda-1$.  

We obtain the path $u^{0,0}_0\,P'\,u^{0,\lambda}_{\alpha}$,  from $u_0$ in $H_{0,0}$ to $u_{\alpha}$ in $H_{0,\lambda}$, where  
$u_{\alpha} = u_0$ for $\lambda$ odd,  $u_{\alpha} = u_1$ for $\lambda$ even, that covers all vertices in $K_0\cup K_{\lambda}$,
together with the vertices in $\cup^{\lambda-1}_{j=1}\,U_j$. The path $u^{0,0}_0\,P'\,u^{0,\lambda}_{\alpha}$ yields a hamilton cycle in $G$ if $\lambda=1$, since $u^{0,0}_0$ and $u^{0,\lambda}_{\alpha}=u^{0,1}_{0}$ are adjacent in $G$. Therefore, the assertion is true if $\lambda=1$.
In the following, we consider $\lambda>1$. In this case, the paths $u^{0,1}_0\,P^*\,u^{0,\lambda-1}_{\alpha}$ and  $u^{0,0}_0\,P'\,u^{0,\lambda}_{\alpha}$ partition the vertices of $G$, and we can connect them
through the edges of type $a$ from $u^{0,0}_0$ to $u^{0,1}_0$,  and 
from $u^{0,\lambda-1}_{\alpha}$ to $u^{0,\lambda}_{\alpha}$; we find a hamilton cycle in $G$.  Notice that the above constructions are valid even in the case where some of the vertices that are used in the description of the paths are the same, i.e., if $v_0=v_k$ and $u_1=u_z$ and/or $v_s=v_h$ and $u_t=u_z$.
\end{proof}

\subsection{Construction for the graphs $B(m; R, S, T)$}

In order to develop the second construction of this section, the following remark is in order.

\begin{remark}\label{rem_even_outer_edges_in_ham_cycle}
    In a connected bicirculant graph, a hamilton cycle  containing a given number of outer edges also contains the same number
    of inner edges. In fact, a hamilton cycle $C$ is the vertex disjoint union of its outer and inner subpaths, which are joined by the spokes to form $C$.
    Let $\mathcal O$ and $\mathcal I$ be the subgraphs of $C$ whose connected components are the outer and inner 
    subpaths of $C$, respectively. Since the end-vertices of the paths in $\mathcal O$ are matched with the end-vertices of the paths in $\mathcal I$ by the spokes, the subgraphs $\mathcal O$ and $\mathcal I$ have the same number of components, say $\omega$. Moreover, $|V(\mathcal O)|=|V(\mathcal I)|$ because
    $C$ has the same number of outer and inner vertices. Since $|V(\mathcal O)|=|E(\mathcal O)|+\omega$, and $|V(\mathcal I)|=|E(\mathcal I)|+\omega$, we have $|E(\mathcal O)|=|E(\mathcal I)|$, that is, $C$ contains the same number of outer and inner edges. 
    
\end{remark}

We now define the second construction.

\begin{lemma}\label{pro_construction_Bd1bis}
Let $G = B(m; R, S, T)$ be a connected bicirculant graph of degree $d \geq 4$, with $|R|,|T|\geq 3$.
Let $H$ be the subgraph of $G$ obtained by removing the edges of types $a\in R$ and $b\in T$, with $a, b\neq m/2$.

If the connected components of $H$ contain a hamilton cycle with at least two outer edges, then $G$ is hamiltonian.
\end{lemma}

\begin{proof}
The assertion is true if $H$ is connected.  Let us assume that $H$ is disconnected, and let $H_i$, with $0\leq i\leq\gcd(m, R\smallsetminus\{a,-a\}, S, T\smallsetminus\{b,-b\})-1$, 
be its connected components; $H_0$ is the component containing $u_0$. Let $C$ be a hamilton cycle in $H_0$ containing at least two outer edges.
We also consider the subgraph $K$ of $G$ obtained by removing the edges of type $a$ from $G$, and denote its connected components with $K_0,\ldots, K_{\lambda}$,  with $\lambda=\gcd(m, R\smallsetminus\{a,-a\}, S, T)-1$; $K_0$ is the component containing $u_0$. We can repeat the same arguments as at the beginning of Section \ref{subsection:construction1}, and say that
the components of $K$ partition the components of $H$: each component of $K$ consists of exactly $\mu+1=\gcd(m, R\smallsetminus\{a,-a\},  S, T\smallsetminus\{b,-b\})/\gcd(m,R\smallsetminus\{a,-a\}, S, T)$ components of $H$,  which are joined by edges of type $b$.

We distinguish two cases: $\mu=0$ and $\mu>0$.

\medskip

\noindent
{\bf Case 1: $\mu =0$. }
Each component of $K$ contains exactly one component of $H$.  We can set $H_j\subseteq K_j$, for $0\leq j\leq\lambda$.
Note that $\lambda>0$ since we are dealing with the case where $H$ is diconnected, i.e., $\gcd(m, R\smallsetminus\{a,-a\},  S, T\smallsetminus\{b,-b\})>1$, and in this case we have $\lambda=\gcd(m, R\smallsetminus\{a,-a\},  S, T\smallsetminus\{b,-b\})-1$.
For vertices and paths in $K_j$, we will use the same notation as at the beginning of Section %\ref{sec:hamconstrucion}.
\ref{subsection:construction1}.
More specifically,  for $0\leq j\leq\lambda$, we denote with $C_{j}$ the cycle obtained from $C$ by adding $ja\pmod m$ to the subscripts of the vertices in $C$; $C_{0}$ corresponds to the cycle $C$. 
The vertices in $C_{j}$ will be denoted with $u^{j}_x$, $v^{j}_x$, where $u^{j}_x=u_{x+ja}$, $v^{j}_x=v_{x+ja}$; the outer vertices $u^{j}_x$ in $C_{j}$ are adjacent to the outer vertices $u^{j+1}_x$ in $C_{j+1}$ by edges of type $a$  (where $C_{\lambda+1} = C_0).$  Moreover,  we will denote with $u_x\,C_{j}\,u_y$ a path corresponding to the subpath $u_x\,C\,u_y$ of $C$ from the vertex $u_x$ to the vertex $u_y$ of $C$, where the vertices $u_x$, $u_y$ in the notation $u_x\,C_{j}\,v_y$ correspond to the vertices $u^{j}_x$, $u^{j}_y$.  Which of the two possible paths we take will be clear from the context. We can thus connect the outer vertices in $u_x\,C_{j}\,u_y$ to the outer vertices in $u_x\,C_{j+1}\,u_y$ that are assigned the same label $u_x$ (or $u_y$) for every $0\le j\le\lambda-1$. We exclude the case with $j=\lambda$ and $j+1=0$ because the vertices $u^{\lambda}_x$, $u^0_x$ might not be adjacent.

Without loss of generality we can assume that $u_0$ is an endvertex of an outer edge that is in $C$. Let $u_0u_s$, $u_hu_{k}$ be distinct outer edges in $C$; their removal  yields two paths, denoted with $u_0\,C\,u_h$ and $u_{s}\,C\,u_{k}$,  that partition the vertices of $C$ (if $u_0=u_h$, then $u_0\,C\,u_h$ consists of a single vertex, the same holds for $u_{s}\,C\,u_{k}$, if $u_s=u_k$).
We also denote with $u_0\,C\,u_s$ and $u_h\,C\,u_k$ the hamilton paths of $H_0$ that are obtained from $C$ by deleting the edges $u_0u_s$ and $u_hu_k$, respectively.

We construct a hamilton cycle in $G$ as follows.  In each $K_j$, with $0< j<\lambda-1$, we take the paths  $u_0\,C_j\,u_h$ and $u_{s}\,C_j\,u_{k}$; we consider 
$u_0\,C_{0}u_s$ in $K_0$, and $u_0\,C_{\lambda}u_s$ or $u_h\,C_{\lambda}u_k$ in $K_{\lambda}$, according to whether $\lambda$ is odd or even, respectively.

We connect the path $u_0\,C_{0}u_s$ to the paths $u_0\,C_1\,u_h$ and $u_{s}\,C_1\,u_{k}$ by adding the edges $u_0^0u_0^1$ and $u_s^0u_s^1$.
For $0< j<\lambda-1$,  $j$ odd (respectively,  $j$ even), we connect the path $u_0\,C_j\,u_h$ to the path $u_0\,C_{j+1}\,u_h$ by joining the vertices that are assigned the label $u_h$ in $C_j$ and $C_{j+1}$
(respectively,  $u_0$); we also connect the path $u_{s}\,C_j\,u_{k}$ to $u_{s}\,C_{j+1}\,u_{k}$ by joining the vertices that are assigned the label $u_{k}$ in $C_j$ and $C_{j+1}$ (respectively,  $u_s$).
We thus find a path $P$ of $G$ from $u^{\lambda-1}_0$ to $u^{\lambda-1}_s$ if $\lambda$ is odd and from $u^{\lambda-1}_h$ to $u^{\lambda-1}_k$ if $\lambda$ is even; path $P$ covers all vertices of $G$, with the exception of the vertices in $K_{\lambda}$.
To construct a hamilton cycle in $G$, all that remains is to connect the end vertices of $P$ to the end vertices of  $u_0\,C_{\lambda}u_s$ or $u_h\,C_{\lambda}u_k$ that are assigned the same label
($u_0$ and $u_s$ for $\lambda$ odd, $u_h$ and $u_k$ for $\lambda$ even).  Figure \ref{fig_construction_Bd1} summarizes the construction for $\lambda=3$.

\begin{figure}[ht]
\begin{center}
\includegraphics[width=12cm]{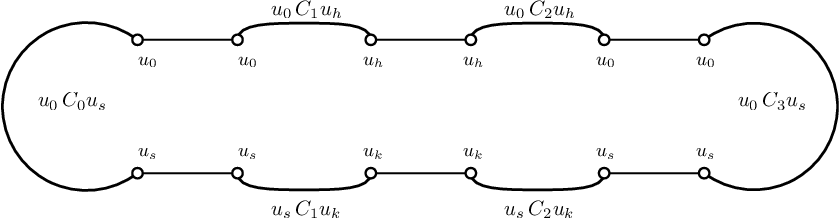}
\end{center}
\caption{The construction of a hamilton cycle described in Lemma  \ref{pro_construction_Bd1bis}; the figure shows the case $\lambda=3$.
The paths we consider in each component $K_j$, $0\leq j\leq\lambda$, are represented by curved lines, the straight line segments represent the edges of type $a$ connecting the vertices in $K_j$ to the vertices in 
$K_{j+1}$ that are assigned the same label.}\label{fig_construction_Bd1}
\end{figure}

\medskip

\noindent
{\bf Case 2: $\mu>0$.}
We show that $K_0$ -- and hence each component $K_j$ with $0< j\leq\lambda$ -- is hamiltonian.
By Remark \ref{rem_even_outer_edges_in_ham_cycle},  the cycle $C$ of $H_0$ contains at least two inner edges.
We can thus apply the construction of Case 1 to $K_0$ by replacing the outer vertices with the inner vertices, and the outer edges of type $a$ with the inner edges of type $b$. We find a hamilton cycle $C'$ in $K_0$.  In turn,  the cycle $C'$ contains at least two outer edges.
We can thus apply the construction of Case 1 to the cycle $C'$, and find a hamilton cycle in $G$ in the case when $\lambda>0$.
When $\lambda=0$, the cycle $C’$ is already a hamilton cycle in the graph $G=K_0$.

\end{proof}

The next lemma will be fundamental to prove the main results. It shows how to combine the constructions in Lemmas  \ref{pro_method1_case1} and  \ref{pro_construction_Bd1bis} to find a hamilton cycle in a 
more general bicirculant graph.

\begin{lemma}\label{lem_combination}
    Let $G=B(m, R,S,T)$,  with $|S|\ge 2$,  be a connected bicirculant graph, and let $H=H(m; S)$ be the cyclic Haar graph induced by the spokes of $G$.
    Assume that $H$ is   disconnected and its connected components are hamiltonian. 
    Further to the above,  assume that $R\cup T$ contains at least one element  that is different from $m/2$ and is not coprime to $\gcd(m,S)$.
    Then $G$ is hamiltonian.
\end{lemma}

\begin{proof}
We can consider $m>5$,  since the bicirculant graphs with $m\le 5$ 
and $|S|\ge 2$ are hamiltionian by Lemma \ref{lemma_bicirculant_m<5}.

Let $L=B(m; R', S,  T')$ be a connected regular spanning subgraph of $G$  that  contains $H$ as a subgraph, 
and the removal of all outer edges of type $a'$ and inner edges of type $b'$ for any pair $a' \in R'$, $b' \in T'$ (with $a',b' \ne m/2$ or $a'=b'=m/2)$ yields a disconnected graph.  Such a graph $L$ always exists: we start from the graph $G$ and remove outer and inner edges of given types so that at each step the graph remains connected until we find a connected spanning subgraph of $G$ with the required property. 
We first remove the outer/inner edges of types $a' \in R$, $b' \in T$ that are coprime to $\gcd(m,S)$. In this way we ensure that at least one element from $R' \cup T'$  that is not coprime to $\gcd(m,S)$ is contained in $R' \cup T'$; without loss of generality we may assume that such an element is from the set $T'$ and we denote it by $b$. 
    
We set $|S|=s$, and denote the degree of $L$ with $d'$. Notice that $d'-s>0$, since we assume $H$ is disconnected, 
and $L$ is defined as a connected graph. We distinguish  the cases $d'-s$ even, and $d'-s$ odd.
    
Assume that $d'-s$ is even. Then there are no edges of type $m/2$ in $R',T'$. We set $d'-s=2k>0$, $R'=\{a_i,-a_i : 1\leq i\leq k\}$ and $T'=\{b_i, -b_i : 1\leq i\leq k\}$, where $b_1=b$, with $b$ not coprime to $\gcd(m,S)$.
For $0\leq j\leq k$, we also set $R_j=\{a_i, -a_i : 1\leq i\leq j\}$ and $T_j=\{b_i, -b_i: 1\leq i\leq j\}$, and denote with $L_j$ the subgraph of $L$ represented as $B(m; R_j, S, T_j)$; in this notation we have $R_0=T_0=\emptyset$, that is,
$L_0=H$, and $L_{k}=L$.

We now construct a hamilton cycle in $G$.  First, we construct a hamilton cycle in each connected component  of $L_1$ by  using the construction from Lemma \ref{pro_method1_case1}. We can do this since $\gcd(m,S,b)>0$.
The construction from  Lemma \ref{pro_method1_case1} provides hamilton cycles containing at least two outer edges (see Figures \ref{fig_path_lambda_positive_mu0}--\ref{fig_ham_cycle_lambda_mu_even_H(m,S)_hamiltonian}).
Therefore, we can apply Lemma  \ref{pro_construction_Bd1bis} to the components of $L_1$, and find a hamilton cycle in each component of $L_2$.
Since also the construction from Lemma \ref{pro_construction_Bd1bis} yields  hamilton cycles containing at least two outer edges (see Figure \ref{fig_construction_Bd1}), we can apply Lemma  \ref{pro_construction_Bd1bis} to the components of $L_2$,  and find a hamilton cycle in each component of $L_3$. By repeatedly applying the construction from Lemma \ref{pro_construction_Bd1bis},  we find a hamilton cycle in $L_k=L$, and therefore in $G$,  since 
$L$ is a connected spanning subgraph of $G$. Thus, the assertion follows for $d'-s$ even.

Assume now  that $d'-s$ is odd. Then there are edges of type $m/2$ in $L$. 
We set $d'-s=2k+1$, $R'=\{a_i, -a_i : 1\leq i\leq k+1\}$, $T'=\{b_i, -b_i : 1\leq i\leq k+1\}$, with $a_{k+1}=b_{k+1}=m/2$, and define $R_j$, $T_j$, $L_j$,  $b_1$ as above; in this case, $L_{k+1}=L$.

We construct a hamilton cycle $C$ in each component of $L_{k}$ by
repeating the same arguments as for the case $d'-s$ even.  By  definition of $L$ and Lemma \ref{lemma_cprod}, 
the graph $L$ is the cartesian product $L=L_k\square K_2$. The existence of a hamilton cycle in $L$, and therefore in $G$,
follows from the fact that $L$ contains the spanning subgraph $C\square K_2$, which is hamiltonian.
\end{proof}

\section{Bicirculants with $s=1$ and $s=2$}\label{subsection:consequences}

In this section we use the constructions presented in  the previous section to complete the classification of hamiltonian bicirculants with $s=1$ and $s=2$.
{We denote by $\B^d_s$ the class of all $d$-valent bicirculants  for which every vertex is adjacent to $s$ spokes.

\begin{proposition}  \label{lem:hamd11}
Every connected bicirculant graph whose vertices are incident to exactly one spoke  is hamiltonian,
except for the complete graph $K_2$ and the generalized Petersen graphs $G(m, 2)$, with $m\equiv 5\pmod 6$, $m\geq 5$.
\end{proposition}

\begin{proof}

The assertion is trivial for $d\le 2$. For $d=3$, it follows from Theorem  \ref{th_cubicbicirculants}. 

We  now consider $d>3$,  and treat separately the cases when $d$ is odd,  and when $d$ is even. 
Notice that a bicirculant from the class $\B(m;d,1)$ does not contain any inner or outer edges of type $m/2$ when  $d$ is odd,  whereas each vertex is incident to such an edge if $d$ is even.

First, we prove that the graphs from the class $\mathcal B^5_1$ are hamiltonian. Then,  we will use induction to prove that every graph from the class $\mathcal B^d_1$, with $d>5$, is also hamiltonian.
We will make use of the fact that for $m \ge 3$, every hamilton cycle in a graph from the class $\B_1^d$ contains at least two outer edges since in such a graph no two spokes are adjacent.

Let $G = B(m; R, \{0\}, T)$ be a connected graph from the class $\mathcal B^5_1$.
Let $H$ be the subgraph obtained from $G$ by removing the edges of types $a$ and $b$, with $a\in R$ and $b\in T$.  The connected components of $H$
are graphs from the class $\mathcal B^3_1$, and we know that they are hamiltonian, except for the Alspach generalized Petersen graphs.

If the components of $H$ are hamiltonian, then the existence of a hamilton cycle in $G$ follows from Lemma \ref{pro_construction_Bd1bis}.

Assume now that the components of $H$ are Alspach generalized Petersen graphs $G(n, 2)$, with $n\equiv 5\pmod 6$, $n\geq 5$. 
By Theorem 4.2 in \cite{AlspachLiu}, we know that every pair of non-adjacent vertices in $G(n, 2)$ is connected by a hamilton path. 

If $H$ is connected, then we find a hamilton cycle in $G$ by adding the edge $u_0u_a$ of $G$ to a hamilton path in $H$ from $u_0$ to $u_a$
(the existence of a hamilton path in $H$ from $u_0$ to $u_a$ is guaranteed by  \cite[Theorem  4.2]{AlspachLiu}, since $u_0$ and $u_a$ are not adjacent in $H$).

If $H$ is disconnected, then we consider the graph $K$ obtained from $G$ by removing the edges of type $a$.  As remarked in the proof of Lemma \ref{pro_construction_Bd1bis}, each
component of $K$ contains $\mu+1$ components of $H$, which are joined by edges of type $b$, where
$\mu+1=\gcd(m, R\smallsetminus\{a,-a\}, T\smallsetminus\{b,-b\})/\gcd(m,R\smallsetminus\{a,-a\},T)$.

If $\mu=0$, then each component of $K$ contains exactly one component of $H$.  We denote the components of $K$ and $H$ with $K_0,\ldots, K_{\lambda}$, 
and $H_0,\ldots, H_{\lambda}$, respectively, where $\lambda=\gcd(m, R\smallsetminus\{a,-a\}, T)-1$, and assume $H_j\subseteq K_j$ for $0\leq j\leq\lambda$.
In the component $K_0$ containing the vertex $u_0$,  there exists a hamilton path $u_0\,P\,u_a$ from $u_0$ to $u_a$ by \cite[Theorem  4.2]{AlspachLiu},  since $H_0\subseteq K_0$.
We denote with $u^{\lambda}_0$, $u^{\lambda}_a$ the copies of the vertices $u_0$, $u_a$ in $K_{\lambda}$, respectively.
The path $u_0\,P\,u_a$ has at least one outer edge, say $u_hu_k$,  since every vertex is incident to exactly one spoke. 
We can thus apply the construction in Case 1 of Lemma \ref{pro_construction_Bd1bis}, with $u_a=u_s$.  
We find a hamilton cycle in $G$ if $\lambda$ is odd, or a hamilton path from $u^{\lambda}_0$ to $u^{\lambda}_a$ if $\lambda$ is even, 
which we turn into a hamilton cycle of $G$, by adding the edge 
$u^{\lambda}_0u^{\lambda}_a$.

The case $\mu>0$ remains to be considered.  We show that the connected component of $K$ containing $u_0$, say $K_0$,  is hamiltonian.  
The component $K_0$ contains $\mu+1$ copies of $H_0$, where $H_{0}$ is the component of $H$ containing $u_0$.  We denote with 
$v^{\mu}_0$, $v^{\mu}_b$  the copies of $v_0$, $v_b$ in the $(\mu+1)$-th copy of $H_0$ in $K_0$.
By \cite[Theorem  4.2]{AlspachLiu}, there exists a hamilton path in $H_0$ from $v_0$ to $v_b$, say $v_0\,P\,v_b$;
it contains at least one inner edge,  say $v_hv_k$. We apply the construction in Case 1 of Lemma \ref{pro_construction_Bd1bis}, by
replacing the outer vertices $u_0$, $u_s$,$u_h$,$u_k$ with the inner vertices $v_0$, $v_b$,$v_h$,$v_k$, respectively,  and the outer edges of type $a$ with the inner edges of type $b$. We  find a hamilton cycle of $K_0$ if $\mu$ is even,  or a hamilton path of $K_0$ 
from $v^{\mu}_0$ to $v^{\mu}_b$ if $\mu$ is odd, which we turn into a hamilton cycle of $K_0$ by adding the edge $v^{\mu}_0v^{\mu}_b$.
In turn, a hamilton cycle of $K_0$ has at least two outer edges, and we can apply the construction in Case 1 of Lemma \ref{pro_construction_Bd1bis}, to construct a hamilton cycle in $G$.  It is thus proved that every graph from the class $\mathcal B^5_1$ is hamiltonian.

We now use induction on $d$ to prove that every graph from the class $\mathcal B^d_1$, with $d$ odd, $d>5$, is hamiltonian.  We assume that the assertion is true for every connected graph
from the class $\mathcal B^{d-2}_1$, with $d-2\geq 5$, $d$ odd.  Let $G= B(m; R, \{0\}, T)$  be a graph from the class $\mathcal B^d_1$. 
Let $H$ be the subgraph obtained from $G$ by removing the edges of types $a$ and $b$,  with $a\in R$ and $b\in T$.  The connected components of $H$, possibly one,  are graphs from the class $\mathcal B^{d-2}_1$, which are hamiltonian by the induction hypothesis.
By Lemma \ref{pro_construction_Bd1bis}, the graph $G$ is also hamiltonian. 

Let us consider the case with $d$ even, $d\geq 4$.  Let $G=B(m; R, \{0\}, T)$ be a graph from the class $\mathcal B^d_1$. As remarked, such graphs have outer or inner edges of type $m/2$ at each vertex, so $m$ must be even.  Let $H$ be the subgraph obtained from $G$ by removing such edges at each vertex.  The graph $H$ has at most two connected components,  which belong to the class $\mathcal B^{d-1}_1$, with $d-1$ odd,  $d-1\geq 3$.  
From the previous results on the class $\mathcal B^{d-1}_1$ with $d-1$ odd,  the connected components of $H$ are hamiltonian if $d>4$.
The connected components of $H$ are also hamiltonian by Theorem \ref{th_cubicbicirculants} when $d=4$ and they are not Alspach generalized Petersen graphs. 
In both cases, the graph $G$ is hamiltonian by Lemma \ref{pro_construction_Bd1bis}.
However, if a connected component  of $H$ is an Alspach generalized Petersen graph, the graph $H$ must consist of two copies of this graph, otherwise $m$ is not even. Recall that by \cite[Theorem  4.2]{AlspachLiu}, every component of the graph $H$ contains a hamilton path.
Therefore, also the graph $G=H^=$  is hamiltonian by Lemma \ref{lemma_cprodham}.
\end{proof}

\begin{proposition}\label{cor_combination}
Every connected bicirculant graph whose vertices are incident to exactly two spokes is hamiltonian.
\end{proposition}
\begin{proof}
Let $G=B(m; R, S, T)$ be a connected bicirculant graph with $|S|=2$. If $m\le 5$, then 
$G$ is hamiltonian by Lemma \ref{lemma_bicirculant_m<5}. Graph $G$ is hamiltonian also when its spokes induce a connected graph, since in this case $H(m;S)$  is a cycle.
We have a hamilton cycle in $G$ even when $G=B(m; m/2, S,m/2)$ and $\gcd(m,S)>1$, since in this case graph $G$ is the cartesian product of  the cycle on $m$ vertices and $K_2$; see Lemma \ref{lemma_cprodham} (iii).

Therefore,  it remains to consider the case where $m>5$,  the spokes induce a disconnected subgraph, and 
the set $R\cup T$ contains at least one element that is different from $m/2$.

If there exists at least one integer $a\in R\cup T \smallsetminus \{m/2\}$ that is not coprime to $\gcd(m,S)$, then $G$ is hamiltonian by Lemma \ref{lem_combination}.

If every element in $R\cup T \smallsetminus \{m/2\}$ is coprime to $\gcd(m,S)$, then $G$ contains a connected generalized rose window graph spanning its vertices. Since generalized rose window graphs are hamiltonian by Theorem \ref{th_rwgraphs}, graph $G$ is hamiltonian. 
\end{proof}

By combining Propositions  \ref{lem:hamd11} and \ref{cor_combination} we obtain Theorem \ref{th_bicirculant_s2}.

\section{Proofs of main theorems}
\label{section:main}

In Section \ref{subsection:consequences}, we gave a proof of Theorem \ref{th_bicirculant_s2}.
In this section we give  proofs of the other  main results
that were already stated in the introduction, namely Theorems \ref {thm:primes}, \ref{th_m/2} and \ref{th_main}. 
The proof of Theorem \ref {thm:primes} relies on the fact that if $m$ is a product of a small number of prime powers, then a bicirculant of order $2m$ contains a connected subgraph for which we know that is hamiltonian (with the known exceptions). In the proofs of Theorems  \ref{th_m/2} and \ref{th_main} we use constructions from Section \ref{sec:hamconstrucion}. We first state an  auxiliary result from \cite{BoPiZi2025}.

\begin{lemma} \label{lemma:primesHaar}
Let $G=H(m;S)$ be a connected cyclic Haar graph with $|S| \ge 4$. If  $m$ is a product of at most three prime powers, then $G$ is hamiltonian. 
\end{lemma}

We now give proofs of Theorems \ref {thm:primes}, \ref{th_m/2} and \ref{th_main}.
\medskip

\noindent
\emph{Proof of Theorem \ref{thm:primes}.} 
Assume that $m$ is a product of at most three prime powers. If $\gcd(m,S)= 1$, then $G$ contains a connected cyclic Haar graph $B(m; \emptyset, S, \emptyset)$ as a subgraph, which is hamiltonian by Theorem \ref{th_cyclic_Haar_s3}  (when $|S|=3$) or by Lemma  \ref{lemma:primesHaar} (when $|S|>3$).

We  now assume that $\gcd(m,S)=\delta>1$. 
If $G$ contains a connected bicirculant whose vertices are incident to two spokes as a subgraph, then it is hamiltonian by Proposition \ref{cor_combination}.
Therefore we assume that this is not the case and we will show that then $m$ needs to be a product of at least four prime powers. 

Let $S=\{0,c_1,\dots, c_{s-1}\}$, where $s=|S| \ge 3$. Since $\gcd(m,S)=\delta >1$, there exists a prime $r$ that divides $\delta$. The graph $G$ is connected, therefore $\gcd(m,R,S,T)=1$ and thus $\gcd(m;R,T)$ is coprime to $r$. However, we have assumed that $G$ does not contain a connected bicirculant whose vertices are incident to two spokes as a subgraph, therefore $\gcd(m, R,T,c_1)>1$. Thus, there exists a prime $p$ that divides $\gcd(m, R,T,c_1)$ and is coprime to $r$. In particular $p \ne r$.  Since $\gcd(m, R,S,T)=1$, there exists $c_i \in S$ that is not divisible by $p$. Again, $\gcd(m;R,T,c_i)>1$, therefore there exists a third prime, say $q$, that divides $\gcd(m, R,T,c_i)$. Thus $m$ is a product of at least three prime powers, namely the powers of the primes $p,q$ and $r$. 

Suppose that $m$ is a product of exactly three prime powers. Since the graph $G$ is connected, there exists $c_j \in S \smallsetminus \{c_i\}$ that is not divisible by $q$ (it may happen that $c_j=c_1$). Since $\gcd(m;R,T,c_j)>1$, it follows that $c_j$ must be divisible by $p$. Now we have elements $c_i,c_j$ from $S$ such that $c_i$ is divisible by $q$ and is coprime to $p$, and $c_j$ is divisible by  $p$  and is coprime to  $q$. But then $c_i-c_j$ is  not divisible by any of $p,q$ and $\gcd(m,R,T)$ is not divisible by $r$. It follows that $\gcd(m, R,T,c_i-c_j)=1$ and $G$ contains a connected  graph $B(m;R,\{0,c_i-c_j\},T)$ as a subgraph, a contradiction. Therefore $m$ is a product of at least four prime powers.

Since the product of the first four   primes is $210$ and
the  product of the first four odd  primes is $1155$ , we now conclude that every bicirculant graph of order $2m$, with  even  $m< 210$ or odd $m<1155$  is hamiltonian, with the exceptions given in Theorem \ref{th_bicirculant_s2}.

\qed

\medskip

\noindent
\emph{Proof of Theorem \ref{th_m/2}.} 
The assertion follows from Theorem \ref{th_generalized_dihedral} if $\gcd(m,S)=1$, since in this case the graph $H(m;S)$ is a 
connected spanning subgraph of $G$ with a hamilton cycle. For $m=4$, the assertion follows from Lemma \ref{lemma_bicirculant_m<5}.
We now assume that $m>4$, $\gcd(m,S)>1$ and consider the subgraph $G'=B(m; m/2,  S,  m/2)$ of $G$, which has $\gcd(m/2, S)\ge 1$ connected components.
Each component of $G'$ has order $m/\gcd(m/2, S)$ that is even. In fact, if 
$m/\gcd(m/2, S)$ is odd, then $\gcd(m/2, S)$ is divisible by the largest power of $2$ dividing $m$, say $2^r$, with $r\ge 1$. It follows that $m/2$ is divisible by $2^r$, a contradiction since $2^{r-1}$ is the largerst power of $2$ dividing $m/2$. Therefore, the components of $G'$ satisfy the assumptions in Theorem \ref{th_generalized_dihedral} and for this reason they are 
hamilton-connected, and consequently they  are also hamiltonian. 
It follows by Lemma \ref{lemma_cprodham} that  $G$ is hamiltonian if
$\gcd(m/2,S)=1$. In the rest of the proof we consider $\gcd(m/2, S)>1$ and treat separately the case where at least one of the integers in $(R\cup T)\smallsetminus\{m/2\}$ is coprime to $\gcd(m/2, S)$, and the case where none of the integers in $(R\cup T)\smallsetminus\{m/2\}$ is coprime to $\gcd(m/2, S)$. 

\medskip

\noindent
{\bf Case 1: at least one of the integers in $(R\cup T)\smallsetminus\{m/2\}$ is coprime to $\gcd(m/2, S)$.}
Without loss of generality we can assume that $a\in R$ is coprime to $\gcd(m/2, S)$. We set $\lambda'=\gcd(m/2, S)-1$ and denote the components of $G'$ with
$H'_0,\ldots, H'_{\lambda'}$, where $H'_0$ is the component containing the vertex $u_0$. 
Since $a$ is coprime to $\gcd(m/2, S)$, the components of $G'$ are cyclically connected through edges of type $a$.
In our notation, the vertices $u_{ia}$ in $H'_i$ are adjacent to the vertices $u_{(i+1)a}$ in $H'_{i+1}$ and, consequently, 
the vertex $u_0$ in $H'_0$ is adjacent to the vertex $u_{-a}$ in $H'_{\lambda'}$.

As noted,  the components of $G'$ are hamilton-connected since Theorem \ref{th_generalized_dihedral} holds.
Thus there exists a hamilton path from $u_0$ to $u_{m/2}$ in $H'_0$, say $u_0\,P'\,u_{m/2}$, and we can also find the hamilton paths 
$u_{\lambda' a}\,P'\,u_{-a}$ and $u_{m/2+\lambda' a}\,P'\,u_{-a}$ in $H_{\lambda'}$ that connect the vertices $u_{\lambda' a}$ and $u_{m/2+\lambda' a}$ to $u_{-a}$  in the case where the vertex $u_{-a}$ is distinct from $u_{\lambda' a}$ and  $u_{m/2+\lambda' a}$, respectively. 
We will combine these paths to construct a hamilton cycle in $G$. More specifically, in each $H_i$, with $0\le i\le\lambda'-1$, we consider the copy of the path $u_0\,P'\,u_{m/2}$, i.e., the path $u_{ia}\,P'\,u_{m/2+ia}$ (the path itself for $i=0$); 
in $H_{\lambda'}$ we consider $u_{\lambda' a}\,P'\,u_{-a}$ or $u_{m/2+\lambda' a}\,P'\,u_{-a}$, according to whether $\lambda'$ is even or odd, respectively. 
We join the paths as follows: for $0\le i\le\lambda'-1$, $i$ even, we add the edge of type $a$ connecting the vertex $u_{m/2+ia}$ in $H_i$ to the vertex $u_{m/2+(i+1)a}$ in $H_{i+1}$;
for $1\le i\le\lambda'-1$, $i$ odd, we add the edge of type $a$ connecting the vertex $u_{ia}$ in $H_i$ to the vertex $u_{(i+1)a}$ in $H_{i+1}$. 
By adding the edge $u_0u_{-a}$, we find a hamilton cycle in $G$, and the assertion follows.
In the case where $\lambda’$ is even and $u_{-a}=u_{\lambda’ a}$, in the above construction we replace the copy of the path $u_{0}\,P'\,u_{m/2}$ 
in $H_{\lambda’-1}$ with the copy of a hamilton path $u_x\,P\,u_{m/2}$ where $u_x$ is an arbitrary vertex different from $u_0,u_{m/2}$, and we take a hamilton path $u_{-a}\,P\,u_{x+\lambda’a}$ in $H_{\lambda’}$. Analogously, in the case where $\lambda’$ is odd and $u_{-a}=u_{m/2+\lambda’a}$,
we replace the copy of the path $u_{0}\,P'\,u_{m/2}$ in $H_{\lambda’-1}$ with the copy of a hamilton path $u_0\,P\,u_x$ where $u_x$ is an arbitrary vertex different from $u_0,u_{m/2}$, and we take a hamilton path $u_{-a}\,P\,u_{x+\lambda’a}$ in $H_{\lambda’}$. In both cases, we can join the copy of the vertex $u_x$ in $H_{\lambda'-1}$ to the copy of $u_x$ in $H_{\lambda'}$, we can leave invariant the other adjacencies and find a hamilton cycle in $G$. 

\medskip

\noindent
{\bf Case 2: none of the integers in $(R\cup T)\smallsetminus\{m/2\}$ is coprime to $\gcd(m/2, S)$}.
Since $G$ is connected, the inequality $\gcd(m/2, a, S)<\gcd(m/2, S)$ holds for at least one element $a\in R\cup T$. 
Without loss of generality we can assume that $a\in R$. Then the subgraph $B(m; \{\pm a, m/2\}, S, \{m/2\})$ of $G$ contains $\gcd(m/2, S)/\gcd(m/2, a, S)>1$ components of $G'$
that are cyclically connected through edges of type $a$. In each component of $B(m; \{\pm a, m/2\}, S, \{m/2\})$ we can repeat the same arguments of Case 1 and find a hamilton cycle, say $C'$, that contains at least two outer edges of type $a$.  This fact allows us to apply the construction from Case 1 of the proof of Lemma \ref{pro_construction_Bd1bis} and find a hamilton cycle in each component of $B(m; \{\pm a, m/2\}, S, \{\pm b, m/2\})$, where $b$ is an arbitrary element in $T$. By setting $L_1=B(m; \{\pm a, m/2\}, S, \{\pm b, m/2\})$, we can also repeat the same arguments as in the first part of the proof of Lemma \ref{lem_combination}. We find a hamilton cycle in $G$, and the assertion follows.

\qed

\medskip

\noindent
\emph{Proof of Theorem \ref{th_main}.} 
The assertion is straightforward if $\gcd(m, S)=1$, while it follows from  Lemma \ref{lemma_bicirculant_m<5} if $m\le 5$. 
In the following we consider $m>5$ and $\gcd(m, S)>1$. By the assumptions, the components of $H(m; S)$ are hamiltonian
and $R\cup T$ contains at least one element, say $b$, that is not coprime to $\gcd(m, S)$. 
If $b\ne m/2$, then the existence of a hamilton cycle in $G$ follows from Lemma \ref{lem_combination}. If $b=m/2$, then $m/2\in R, T$ and
the assertion follows from Theorem \ref{th_m/2}.
\qed

\section{Concluding remarks}
\label{sec_opencases}

In this paper we studied hamiltonicity of bicirculant graphs. We presented constructions which showed how to combine hamilton cycles in the connected components of certain subgraphs of a connected bicirculant to a hamilton cycle of the whole graph. 
As a consequence we showed that many classes of bicirculants are hamiltonian and that all bicirculants of order $2m$ are all hamiltonian when $m<210$, % and even or $m<1155$ and odd, 
with known exceptions. The following result shows where to search for the smallest open cases, that is, bicirculant graphs for which the existence of a hamilton cycle is not known.

\begin{corollary}\label{cor_s=3}
Let $G=B(m; R, S, T)$ be a connected bicirculant graph with $m\ge 3$,  $|S|=3$ and such that $R\cup T$ contains at least one element that is not coprime to $\gcd(m, S)$ in the case where $\gcd(m, S)>1$. Then  $G$ is hamiltonian.  
\end{corollary}

\begin{proof}
By Theorem \ref{th_cyclic_Haar_s3}, the components of $H(m; S)$ are hamiltonian.
Then the assumptions in Theorem \ref{th_main} are satisfied and the assertion follows.
\end{proof}

The first open cases are thus bicirculant graphs $B(m; R, S, T)$ with $|S|=3$, $\gcd(m, S)>1$ and all elements in $R\cup T$ being coprime to $\gcd(m, S)$.  Moreover, $m$ has to be product of at least four prime powers and $m/2\not\in R\cup T$. 
The graphs should also not contain a connected generalized rose window graph or a cubic Cayley graph on a dihedral group as a subgraph, since these graphs are all hamiltonian. One of the first cases with even $m$ that is not covered by our results is, for example, the bicirculant  $B(210;30,\{0,14,35\},60)$. It can be checked by computer that this graph is hamiltonian. A smallest open case with odd $m$ is the graph  
$B(1155;105,\{0,33,110\},315)$,  however, the size of this graph already makes searching for a hamilton cycle challenging.

Other `prominent' open cases are cyclic Haar graphs $H(m; S)$ with odd $m$ that is a product of at least four prime powers, 
and $|S|\ge 4$. We consider cyclic Haar graphs as prominent open cases because by Theorem \ref{th_main} we have that the existence of a hamilton cycle in cyclic Haar graphs of degree $d\ge 4$ implies the existence of a hamilton cycle in  every bicirculant graph $B(m; R, S, T)$ with $R\cup T$ containing at least one element that is not coprime to $\gcd(m, S)$ in the case where $\gcd(m, S)>1$.

We also emphasize that it seems possible to construct a hamilton cycle in bicirculant graphs independently of the existence of a hamilton cycle in cyclic Haar graphs of degree $d\ge 4$. However, the construction still seems to leave open cases, particularly when $m$ is odd.

\section*{Acknowledgements}
Simona Bonvicini is a member of GNSAGA of Istituto Nazionale di Alta Matematica (INdAM).
Toma\v{z} Pisanski is supported in part by the Slovenian Research Agency (research program P1-0294 and research projects J1-4351, J5-4596 and BI-HR/23-24-012).
Arjana \v Zitnik is supported in part by the Slovenian Research Agency (research program P1-0294 and research projects J1-3002 and J1-4351).

\bigskip
\noindent
Simona Bonvicini\\
Dipartimento di Scienze Fisiche, Informatiche e Matematiche\\
Universit\`a di Modena e Reggio Emilia, \\
via Campi 213/b, 41126 Modena, Italy\\
email: \texttt{simona.bonvicini@unimore.it}
\\ %[-3mm]

\noindent
Toma\v z Pisanski  \\ 
University of Primorska, FAMNIT \\
Glagolja\v ska 8, 6000 Koper, Slovenia\\
and\\
Institute of Mathematics, Physics, and Mechanics\\
Jadranska 19, 1000 Ljubljana, Slovenia\\
email: \texttt{tomaz.pisanski@upr.si} 
\\ %[-3mm]

\noindent
Arjana \v Zitnik\\ 
University of Ljubljana, Faculty of Mathematics and Physics\\
Jadranska 19, 1000 Ljubljana, Slovenia\\ 
and\\
Institute of Mathematics, Physics, and Mechanics\\
Jadranska 19, 1000 Ljubljana, Slovenia\\
email: \texttt{arjana.zitnik@fmf.uni-lj.si}

\end{document}